\documentclass[12pt]{amsart}
\usepackage{amssymb}
\usepackage{palatino}
\input amssym.def
\usepackage{amsmath, amsfonts,hyperref}
\usepackage{amscd}
\usepackage[mathscr]{eucal}
\newcommand{\Z}{{\mathbb Z}}
\newcommand{\Q}{{\mathbb Q}}
\newtheorem{thm}{Theorem}
\newtheorem{lem}{Lemma}
\newtheorem{cor}{Corollary}
\newtheorem{prop}{Proposition}
\newtheorem{rmk}{Remark}
\newtheorem{defn}{Definition}
\newcommand{\thmref}[1]{Theorem~\ref{#1}}
\newcommand{\propref}[1]{Proposition~\ref{#1}}
\newcommand{\lemref}[1]{Lemma~\ref{#1}}
\newcommand{\corref}[1]{Corollary~\ref{#1}}

\begin{document}

\title[Error term in a Parseval type formula]
{On the error term in a Parseval type formula in the theory of
Ramanujan expansions}

\author{M. Ram Murty and Biswajyoti Saha}

\address{M. Ram Murty\\ $\phantom{mmmmmmmmmm
mmmmmmmmmmmmmmmmmmmmmmmmmmmmmmmmmmm}$
Department of Mathematics, Queen's University,
Kingston, Ontario, K7L 3N6, Canada}
\email{murty@mast.queensu.ca}

\address{Biswajyoti Saha\\ $\phantom{mmmmmmmmmm
mmmmmmmmmmmmmmmmmmmmmmmmmmmmmmmmmmm}$
Institute of Mathematical Sciences, C.I.T. Campus, Taramani, 
Chennai, 600 113, India}
\email{biswajyoti@imsc.res.in}

\subjclass[2010]{11A25,11K65,11N37}

\keywords{average order, Ramanujan expansions,
Parseval type formula, error terms}

\begin{abstract}
 Given two arithmetical functions $f,g$ we derive, under suitable conditions,
 asymptotic formulas with error term, for the convolution sums $\sum_{n \le N}
 f(n) g(n+h)$, building on an earlier work of Gadiyar, Murty and Padma. A key role in our method is played by the theory of Ramanujan
 expansions for arithmetical functions.
\end{abstract}

\maketitle

\section{Introduction and Statement of Theorems}
The past century of mathematics, in particular number theory, has witnessed
a number of developments in many different directions, originating from
different articles by Srinivas Ramanujan. One of these is the theory of
Ramanujan expansions. In 1918, Ramanujan \cite{SR} introduced certain
sums of roots of unity. To be precise, for positive integers $r,n$, he
defined the following sum,
\begin{equation}\label{def}
 c_r(n):=\sum_{a\in (\Z/r\Z)^*}\zeta_r^{an},
\end{equation}
where $\zeta_r$ denotes a primitive $r$-th root of unity.
These sums are now known as Ramanujan sums. Among other significant
properties of Ramanujan sums, we list a few, which can be obtained
from elementary observations. For a more elaborate account on Ramanujan
sums, we refer the reader to the texts \cite{RS,SS} and the survey articles
\cite{LGL,RM2}. We know:
\begin{itemize}

 \item[a)] For any $r,n$, $c_r(n) \in \Z$. This can be
 seen by reading the sum in \eqref{def} as the trace of the
 algebraic integer $\zeta_r^n$.
 
 \item[b)] For fixed $n$, $c_r(n)$ is a multiplicative function i.e.
 for $r_1,r_2$ with $\gcd (r_1,r_2)=1$ we have
 $c_{r_1r_2}(n)=c_{r_1}(n) c_{r_2}(n)$. This is essentially due to the
 fact that, for $r_1,r_2$ with $\gcd (r_1,r_2)=1$, the fields
 $\Q(\zeta_{r_1})$ and $\Q(\zeta_{r_2})$ are linearly disjoint.
 
 \item[c)] $c_r(\cdot)$ is a periodic function with period $r$.
 In fact, $c_r(n)=c_r(\gcd (n,r))$.
 
 \item[d)] $c_r(n)$ can be expressed in terms of the M\"obius function and
 written as $$c_r(n)=\sum_{d|\gcd(n,r)} \mu(r/d) d.$$
 
\end{itemize}
Ramanujan used these sums to derive point-wise convergent series expansion
of various arithmetical functions of the form $\sum_r a_rc_r(n)$,
which are now called Ramanujan expansions. More precisely, given an
arithmetical function $f$, we say $f$ admits a Ramanujan expansion, if
$$f(n)= \sum_r \hat f(r)c_r(n)$$ for appropriate complex numbers $\hat f(r)$
and the series on the right hand side converges.
Existence of such expansions for a given arithmetical function and their
convergence properties have been studied extensively in the past, for example in
\cite{JS,AH,SS}. However, we do not
discuss these here. In this article we focus on a different theme.

In \cite{GMP}, Gadiyar, Murty and Padma have studied sums of the kind
$\sum_{n \le N} f(n) g(n+h)$ for two arithmetical functions $f,g$ with
absolutely convergent Ramanujan expansions. They derived asymptotic
formulas which are analogous to Parseval's formula in the case of Fourier
series expansions. However it seems that the study of error term for
these sums has not been carried out before. Under certain additional
hypotheses we extend their results and provide explicit error terms.
To be precise we prove,

\begin{thm}\label{error1}
 Suppose that $f$ and $g$ are two arithmetical functions with absolutely
 convergent Ramanujan expansion:
 $$f (n) = \sum_r \hat f(r) c_r(n), \phantom{mm} g(n) = \sum_s \hat g(s)c_s(n),$$
 respectively. Further suppose that
 $$\big|\hat f(r)\big|,\big|\hat g(r)\big| \ll \frac{1}{r^{1 +\delta}}$$
 for some $\delta > 1/2$. Then, we have,
 $$\sum_{n \le N} f(n)g(n) = N \sum_r \hat f(r) \hat g(r) \phi(r) +
 O\left( N^{\frac{2}{1+2\delta}} (\log N)^{\frac{5+2\delta}{1+2\delta}}\right).$$
\end{thm}

\begin{thm}\label{error2}
 Let $f$ and $g$ be two arithmetical functions with the same hypotheses as
 in \thmref{error1} and $h$ be a positive integer. Then we have,
 $$\sum_{n \le N} f(n)g(n+h) = N \sum_r \hat f(r) \hat g(r) c_r(h) +
 O\left( N^{\frac{2}{1+2\delta}} (\log N)^{\frac{5+2\delta}{1+2\delta}}\right).$$
\end{thm}

In his article \cite{SR}, Ramanujan showed that, for real variable $s > 0$,
\begin{equation}\label{sigma_s}
 \frac{\sigma_s(n)}{n^s}=\zeta(s+1)\sum_r \frac{c_r(n)}{r^{s+1}},
\end{equation}
where $\sigma_s(n)=\sum_{d|n}d^s$. Hence as an immediate corollary to
\thmref{error2} we get,

\begin{cor}\label{c1}
 For $s,t > 1/2$ and any non-negative integer $h$, we have,
 $$\sum_{n \le N} \frac{\sigma_s(n)}{n^s}~\frac{\sigma_t(n+h)}{(n+h)^t}
 =N\frac{\zeta(s+1)\zeta(t+1)}{\zeta(s+t+2)} \sigma_{-(s+t+1)}(h)
 +O\left( N^{\frac{2}{1+2\delta}} (\log N)^{\frac{5+2\delta}{1+2\delta}}\right),$$
 where $\delta=\min \{s,t\}$.
\end{cor}
This result, without the explicit error term was mentioned by Ingham in his article
\cite{AEI}. A proof of Ingham's result was obtained in \cite{GMP}.

Ramanujan \cite{SR} also showed that, for $s>0$
\begin{equation}\label{phi_s}
 \frac{\phi_s(n)}{n^s} \zeta(s+1)= \sum_r\frac{\mu(r)}{\phi_{s+1}(r)}c_r(n),
\end{equation}
where $\phi_s(n):= n^s \prod\limits_{\substack{p|n\\ p~prime}}(1-p^{-s})$. Note that,
$\frac{1}{\phi_{s}(r)} \ll \frac{1}{r^s}$. Hence as a corollary we have,

\begin{cor}\label{c2}
 If $s,t \ge 1/2$ and $h$ is a non-negative integer, then
 $$\sum_{n \le N} \frac{\phi_s(n)}{n^s} \frac{\phi_t(n+h)}{(n+h)^t}
 =N \Delta(h) +O\left( N^{\frac{2}{1+2\delta}}
 (\log N)^{\frac{5+2\delta}{1+2\delta}}\right),$$
 where $\Delta(h)=\prod_{p|h} \left[\left(1-\frac{1}{p^{s+1}}\right)\left(1-\frac{1}{p^{t+1}}\right)
 + \frac{p-1}{p^{s+t+2}}\right]
 \prod_{p\nmid h}\left[\left(1-\frac{1}{p^{s+1}}\right)\left(1-\frac{1}{p^{s+1}}\right)
 -\frac{1}{p^{s+t+2}}\right]$ and $\delta=\min \{s,t\}$.
\end{cor}

\section{Preliminaries}
For the sake of completeness we collect a few basic results about certain
arithmetical functions. They can be found easily in the texts like \cite{TA,RM1,SS}.
We first list the results about the average order of  Euler's phi function
$\phi(\cdot)$ and the M\"obius function $\mu(\cdot)$.

\begin{prop}\label{phi}
 For any real number $x \ge 1$,
 $$\sum_{k \le x} \phi(k)= \frac{3}{\pi^2}x^2 + O(x \log x).$$
\end{prop}

\begin{defn}
 \rm
 The Mertens function $M(\cdot)$ is defined for all positive integers $n$ as
 $$M(n) := \sum_{k\le n} \mu(k)$$
 where $\mu(\cdot)$ is the M\"obius function. The above definition can be extended
 to any real number $x \ge 1$ by defining,
 $$M(x) := \sum_{k \le x} \mu(k).$$ 
\end{defn}

Essentially, the error term in the prime number theorem, due to
de la Vall\'ee-Poussin \cite{VP} gives us,

\begin{prop}\label{mertens}
 For any real number $x \ge 1$,
 $$ M(x)= \sum_{k \le x} \mu(k) = O\left(x e^{-c \sqrt{\log x}}\right),$$
 where $c$ is some positive constant.
\end{prop}

\begin{rmk}
 \rm
 This error term has been improved independently by Korobov \cite{NMK}
 and Vinogradov \cite{IMV} in 1958 to
 $O\bigg(x e^{-c \frac{(\log x)^{3/5}}{(\log \log x)^{1/5}}}\bigg)$ with $c >0$,
 the best known to date. An equivalent statement of the Riemann hypothesis is
 $M(x)=O\left(x^{\frac{1}{2}+ \epsilon}\right)$ for any $\epsilon >0$.
\end{rmk}

Next we note the following very useful theorem known as {\it ``partial summation''}.

\begin{thm}\label{partial}
 Suppose $\{a_n\}_{n=1}^\infty$ is a sequence of complex numbers and $f(t)$
 is a continuously differentiable function on $[1, x]$. Let $A(t):=\sum_{n \le t} a_n$.
 Then, $$ \sum_{n \le x} a_n f(n)= A(x)f(x)- \int_1^x A(t)f'(t) ~ dt.$$
\end{thm}

Let $d_k(n)$ be the number of ways writing $n$ as a product of $k$ numbers. Sometimes we write $d(n)$ to denote $d_2(n)$.
Using \thmref{partial} we can derive the following result about the average
order of the arithmetical function $d_k(\cdot)$.

\begin{prop}\label{k-product}
 $$\sum_{n \le x}d_k(n)= \frac{x (\log x)^{k-1}}{(k-1)!}+ O\left(x (\log x)^{k-2}\right).$$
\end{prop}

\bigskip

We will make use of all these results from basic analytic number theory
to extend the following theorems of \cite{GMP}.

\begin{thm}[Gadiyar-Murty-Padma]\label{w/oerror1}
 Suppose that $f$ and $g$ are two arithmetical functions with absolutely
 convergent Ramanujan expansion:
 $$f (n) = \sum_r \hat f(r) c_r(n), \phantom{mm} g(n) = \sum_s \hat g(s)c_s(n),$$
 respectively. Suppose that
 $$\sum_{r,s} \big| \hat f(r) \hat g(s)\big|\gcd(r,s)d(r)d(s) < \infty.$$ Then,
 as $N$ tends to infinity,
 $$\sum_{n \le N} f(n)g(n) \sim N \sum_r \hat f(r) \hat g(r) \phi(r).$$
\end{thm}

\begin{thm}[Gadiyar-Murty-Padma]\label{w/oerror2}
 Suppose that $f$ and $g$ are two arithmetical functions as in
 \thmref{w/oerror1} and $h$ is a positive integer. Suppose further that
 $$\sum_{r,s} \big| \hat f(r) \hat g(s)\big|(rs)^{1/2}d(r)d(s) < \infty .$$
 Then, as $N$ tends to infinity,
 $$\sum_{n \le N} f(n)g(n+h) \sim N \sum_r \hat f(r) \hat g(r) c_r(h).$$
\end{thm}

\begin{rmk}
\rm Our hypotheses about $\hat f ,\hat g$,
 $$\big|\hat f(r)\big|,\big|\hat g(r)\big| \ll \frac{1}{r^{1 +\delta}}
 ~\mbox{ for }~ \delta > 1/2,$$ are extensions of the condition
 $$\sum_{r,s} \big| \hat f(r) \hat g(s)\big|(rs)^{1/2}d(r)d(s) < \infty,$$
 and include it as a consequence. Also note that, if we want
 to extend the hypothesis of \thmref{w/oerror1}, \ref{w/oerror2}
 in the form that we have in \thmref{error1}, \ref{error2} then
 $\delta > 1/2$ is the optimal choice.
\end{rmk}

To prove these theorems they prove certain lemmas about the sums of the from
$$\sum_{n \le N} c_r(n) c_s(n+h).$$ We will also make use the following.

\begin{lem}\label{l1}
 $$\sum_{n \le N} c_r(n)c_s(n+h)=\delta_{r,s} N c_r(h) + O(rs \log rs),$$
 where $\delta_{\cdot,\cdot}$ denotes the Kronecker delta function.
\end{lem}

\begin{lem}\label{l2}
 $$\Big |\sum_{n \le N} c_r(n)c_s(n+h)\Big |\le d(r) d(s) \sqrt{rs N(N+h)}.$$
\end{lem}

There are other related results in \cite{GMP} which are of independent interest
but not relevant to our work here.

\section{Proofs of the Theorems}

\subsection{Proof of \thmref{error1}} We start as it is done in \cite{GMP}.
Let $U$ be a parameter tending to infinity which is to be chosen later. We have
by absolute convergence of the series,
\begin{align*}
 \sum_{n \le N} f(n) g(n) &= \sum_{n \le N} \sum_{r,s} \hat f(r) \hat g(s) c_r(n)c_s(n)\\
 &= A+B, ~\mbox{ where }
\end{align*}
$$ A= \sum_{n \le N} \sum_{\substack{r,s\\ rs \le U}} \hat f(r) \hat g(s) c_r(n)c_s(n)
~\mbox{ and }~
B=\sum_{n \le N} \sum_{\substack{r,s\\ rs > U}} \hat f(r) \hat g(s) c_r(n)c_s(n).$$
Interchanging summations and applying \lemref{l1} (for $h=0$) we get,
\begin{align*}
 A&= N \sum_{r^2 \le U} \hat f(r) \hat g(r) \phi(r) + O(U \log U)\\
 &=C+D+O(U \log U), ~\mbox{ where }
\end{align*}
$$C=N \sum_r \hat f(r) \hat g(r) \phi(r) ~\mbox{ and }~
D=-N \sum_{r^2 > U} \hat f(r) \hat g(r) \phi(r).$$
Note that $C$ is the main term according to our theorem. Using the hypothesis we get,
$$D = O \left( N \sum_{r > \sqrt U}\frac{\phi(r)}{r^{2+2\delta}}\right).$$
Using \thmref{partial} and \propref{phi} we get,
$$O\left(\sum_{r > \sqrt U}\frac{\phi(r)}{r^{2+2\delta}}\right)
=O\left(\frac{1}{U^\delta} + \int_{\sqrt U}^\infty \frac{t^2}{t^{3+2\delta}}~dt\right)
=O\left( \frac{1}{U^\delta} \right).$$
Hence we obtain, $D=O \left( \frac{N}{U^\delta}\right)$.

\bigskip

Now, for $B$, interchanging the summation and applying \lemref{l2} (for $h=0$) we get,
\begin{align*}
 B& \ll N \sum_{rs > U} \frac{d(r) d(s) \sqrt{rs}}{(rs)^{1+\delta}}\\
 &=N \sum_{rs > U} \frac{d(r) d(s)}{(rs)^{1+(\delta - 1/2)}}
\end{align*}
Note that, $\sum_{rs=t}d(r)d(s)=d_4(t)$. Using \thmref{partial} and \propref{k-product} we get,
\begin{align*}
 O\left( \sum_{t > U} \frac{d_4(t)}{t^{1+(\delta - 1/2)}} \right )
 &= O\left(\frac{U(\log U)^3}{U^{1+(\delta-1/2)}} + \int_U^\infty
 \frac{t(\log t)^3}{t^{2+(\delta-1/2)}}~dt\right )\\
 &= O\left(\frac{(\log U)^3}{U^{(\delta-1/2)}}+ \int_U^\infty
 \frac{(\log t)^2}{t^{1+(\delta-1/2)}}~dt\right ), ~\mbox{ integrating by parts }\\
 &= O\left(\frac{(\log U)^3}{U^{(\delta-1/2)}}\right ),
 ~\mbox{ integrating by parts multiple times.}\\
\end{align*}
Hence we end up getting, $B=O\left(\frac{N(\log U)^3}{U^{(\delta-1/2)}}\right )$ and thus
$$\sum_{n \le N} f(n) g(n)=N \sum_r \hat f(r) \hat g(r) \phi(r) + O(U \log U)
+O\left(\frac{N(\log U)^3}{U^{(\delta-1/2)}}\right ).$$
To optimize the error term we choose,
$U=N^{\frac{2}{1+2\delta}} (\log N)^{\frac{4}{1+2\delta}}$
and for this choice of $U$ we get
$$\sum_{n \le N} f(n) g(n)=N \sum_r \hat f(r) \hat g(r) \phi(r)
+ O\left( N^{\frac{2}{1+2\delta}} (\log N)^{\frac{5+2\delta}{1+2\delta}}\right).$$
This concludes the proof of \thmref{error1}.

\subsection{Proof of \thmref{error2}} This proof also starts off similarly and we get
$$\sum_{n \le N} f(n) g(n+h) = A+B, ~\mbox{ where }$$
$$ A= \sum_{n \le N} \sum_{\substack{r,s\\ rs \le U}} \hat f(r) \hat g(s) c_r(n)c_s(n+h)
~\mbox{ and }~
B=\sum_{n \le N} \sum_{\substack{r,s\\ rs > U}} \hat f(r) \hat g(s) c_r(n)c_s(n+h).$$
Likewise, the interchange of the summations and \lemref{l1} yield,
$A=C+D+O(U \log U)$, where $C=N \sum_r \hat f(r) \hat g(r) c_r(h)$, the main term and
$D=-N \sum_{r^2 > U} \hat f(r) \hat g(r) c_r(h)$. This time we have,
$$D = O \left( N \sum_{r > \sqrt U}\frac{c_r(h)}{r^{2+2\delta}}\right).$$
To apply \thmref{partial} we need to know about $\sum_{r \le x} c_r(h)$. We write,
$$\sum_{r \le x} c_r(h)=\sum_{r \le x} \sum_{d|r,d|h}\mu(r/d)d
=\sum_{\substack{k,d\\dk \le x, d| h}}d \mu(k)
=\sum_{d|h}d\sum_{k \le x/d} \mu(k).$$
The innermost sum is $M(x/d)$. Now due to \propref{mertens} we get,
$$\sum_{r \le x} c_r(h)=O\left(\sum_{d|h} x e^{-c\sqrt{\log (x/d)}} \right)
=O\left( x e^{-c\sqrt{\log x}}\epsilon(h)\right),$$
for some function $\epsilon(\cdot)$  of $h$ which is bounded above by
$e^{c\sqrt{\log h}}d(h)$. Hence using \thmref{partial} we obtain,
\begin{align*}
 D&=O\left( N \epsilon(h) \left[\frac{\sqrt U e^{-c\sqrt{\log \sqrt U}}}{U^{1+\delta}}
 +\int_{\sqrt U}^\infty\frac{t e^{-c\sqrt{\log t}}}{t^{3+2\delta}}~dt\right]\right)\\
 &=O\left( N \epsilon(h) \left[\frac{1}{U^{1/2+\delta}}
 +\int_{\sqrt U}^\infty\frac{1}{t^{2+2\delta}}~dt\right]\right)
 =O\left(\frac{N \epsilon(h)}{U^{1/2+\delta}}\right).
\end{align*}
For $B$, we apply \lemref{l2}. A similar calculation yields,
$B=O\left(\frac{\sqrt{N(N+h)}(\log U)^3}{U^{(\delta-1/2)}}\right )$.
Hence for fixed $h$ we can write,
$$\sum_{n \le N} f(n) g(n+h)=N \sum_r \hat f(r) \hat g(r) c_r(h) + O(U \log U)
+O\left(\frac{N(\log U)^3}{U^{(\delta-1/2)}}\right )$$
and then similarly as before, choosing
$U=N^{\frac{2}{1+2\delta}} (\log N)^{\frac{4}{1+2\delta}}$ we conclude,
$$\sum_{n \le N} f(n)g(n+h) = N \sum_r \hat f(r) \hat g(r) c_r(h) +
O\left( N^{\frac{2}{1+2\delta}} (\log N)^{\frac{5+2\delta}{1+2\delta}}\right).$$

\subsection{Proof of \corref{c1}} Note that using \eqref{sigma_s} we get,
$$\sum_r \frac{c_r(h)}{r^{s+t+2}}=\frac{1}{\zeta(s+t+2)}
\frac{\sigma_{s+t+1}(h)}{h^{s+t+1}}=\frac{\sigma_{-(s+t+1)}(h)}{\zeta(s+t+2)}.$$
This completes the proof of \corref{c1}.

\subsection{Proof of \corref{c2}} Since $\mu(r),\phi_s(r),c_r(h)$ are multiplicative
functions of $r$ and the M\"obius function is supported at the square free numbers,
we get,
$$\sum_r\frac{\mu^2(r)}{\phi_{s+1}(r)\phi_{t+1}(r)}c_r(h)
=\prod_{p~prime}\left(1+ \frac{\mu^2(p)}{\phi_{s+1}(p)\phi_{t+1}(p)}c_p(h)\right).$$
Now, $c_p(h)=\begin{cases}
              p-1 &\mbox{ if } p|h,\\
              -1 &\mbox{ if } p\nmid h.
             \end{cases}$
Hence we obtain,
\begin{align*}
 &\frac{1}{\zeta(s+1)\zeta(t+1)}\sum_r\frac{\mu^2(r)}{\phi_{s+1}(r)\phi_{t+1}(r)}c_r(h)\\
 =&\prod_{p|h} \left[\left(1-\frac{1}{p^{s+1}}\right)\left(1-\frac{1}{p^{t+1}}\right)
 + \frac{p-1}{p^{s+t+2}}\right]\\
 \times&\prod_{p\nmid h}\left[\left(1-\frac{1}{p^{s+1}}\right)\left(1-\frac{1}{p^{s+1}}\right)
 -\frac{1}{p^{s+t+2}}\right].
\end{align*}
This completes the proof of \corref{c2}.

\bigskip

\noindent{\bf Acknowledgments:} We thank Sanoli Gun and Purusottam Rath for their
comments on an earlier version of this article.

\end{document}